\documentclass[11pt]{article}

\usepackage[top=50mm, bottom=50mm, left=50mm, right=50mm]{geometry}
\usepackage{lineno}
\usepackage{natbib}
\usepackage{amssymb}
\usepackage{amsmath}
\usepackage{amsthm}
\usepackage{epsfig}
\usepackage{graphicx}
\usepackage{graphics}
\usepackage{float}
\usepackage{multirow}
\usepackage{color}
\usepackage{lineno}
\usepackage{fullpage}
\usepackage[normalem]{ulem} 
\usepackage{makeidx}
\usepackage{xspace}
\usepackage{wrapfig}
\usepackage{hyperref}
\usepackage{tikz}
\usetikzlibrary{arrows.meta,positioning,calc, decorations.pathmorphing, decorations.pathreplacing, decorations.shapes}
\usepackage{bbm}
\usepackage{subcaption}
\makeindex

\newtheorem{theorem}{Theorem}

\definecolor{darkred}{rgb}{1, 0.1, 0.3}
\definecolor{darkblue}{rgb}{0.1, 0.1, 1}
\definecolor{darkgreen}{rgb}{0,0.6,0.5}

\newcommand{\prob}{\mathbb{P}}

\title{Epidemics with avoidance and isolation on $\mathbb{Z}^d$}
\author{
{Andrew Heeszel\footnote{Department of Mathematics, University of California, Davis}} \and {Matthew Wascher\footnote{Department of Mathematics, Applied Mathematics, and Statistics, Case Western Reserve University}}
}
\date{}

\begin{document}

\maketitle

\begin{abstract}
The contact process with avoidance is a generalization of the classical contact process (SIS epidemic) that introduces a mechanism for healthy individuals to avoid their infected neighbors. Let $G$ be a directed graph. At each time $t$, each vertex is either healthy or infected and each edge is either active or inactive. Each infected vertex infects each healthy neighbor across each active edge at rate $\lambda$ and recovers at rate $1$, while each active edge pointing from an infected vertex to a healthy vertex becomes inactive at rate $\alpha$. An inactive edge becomes active when its tail vertex recovers. This model has been previously studied on $\mathbb{Z}$, the $n$-cycle $\mathbb{Z}_n$, and the $n$-star graph; here we extend the study of this model to lattices $\mathbb{Z}^d$, $d \geq 2$. We show that for every $d \geq 2$ and fixed $\alpha > 0$, there exist constants $\lambda(\alpha,d)^-$ and $\lambda(\alpha,d)^+$ such that for all $\lambda < \lambda(\alpha,d)^-$ the infection dies out almost surely while for all $\lambda > \lambda(\alpha,d)^+$ the infection persists indefinitely with positive probability. Furthermore, we show that both $\lambda(\alpha,d)^-$ and $\lambda(\alpha,d)^+$ scale like $1/d$ as $d \rightarrow \infty$ and that there exists a constant $C(\alpha)$ such that when $\lambda > C(\alpha)/d$ the process has a nontrivial invariant measure for $d$ sufficiently large. Our methods and most of our results also apply to the SIRS model and a related model in which infected vertices enter an isolated state at rate $\alpha$ and transition from both isolated and infected to healthy at rate $1$.
\end{abstract}

\section{Introduction}\label{sec:intro}

The contact process or SIS epidemic, first introduced in \cite{Harris1974}, is a widely-used model for the spread of a disease through a population. Let $G = (V,E)$ be a graph with vertices $V$ and edges $E$. In the \textit{classical contact process} the state of each vertex $x$ at time $t$ is either $0$ or $1$ where $0$ denotes healthy or susceptible (also denoted by S) and $1$ denotes infected (also denoted by I.) The infection then spreads across each infected-healthy or IS edge according to a rate parameter $\lambda$.

Recently, there has been increased interest in variants of the contact process that include additional vertex or edge states or different transition rules intended to model more complex environmental, biological, or behavioral factors that might affect a real-world epidemic. The model that we study here is the \textbf{contact process with avoidance} introduced in \cite{CSW}. In this model, we take the graph $G = (V,E)$ to have directed edges and include a second parameter $\alpha$ that controls the rate at which each healthy vertex avoids each infected neighbor by temporarily deactivating the IS edge from the infected neighbor to itself. The edge reactivates when the offending vertex recovers.

One important question about the contact process and related models is whether there exist different parameter regimes (or phases) under which the epidemic exhibits qualitatively different behavior. In particular, this question often focuses on the survival time of the epidemic. On infinite graphs, we consider two potential phases. Let $\eta_t$ be the infected set of the classical contact process at time $t$, let $|\eta_t|$ denote the size of this set, let the initial infected set $\eta_0$ be finite, and let $\prob_{\lambda}()$ denote probability given the infection parameter is $\lambda$. The \textit{subcritical} or \textit{extinction phase} is the set of $\lambda$ such that $\prob_{\lambda}(|\eta_t| > 0$ $\forall t) = 0$ while the \textit{supercritical} or \textit{survival phase} is the set of $\lambda$ such that $\prob_{\lambda}(|\eta_t| > 0$ $\forall t) > 0$. Clearly the extinction phase contains $\lambda = 0$ while the survival phase contains $\lambda = \infty$ if one allows this value. We say there is a \textbf{phase transition} if the extinction phase contains some $\lambda \neq 0$ and the survival phase contains some $\lambda \neq \infty$.

While the condition that $|\eta_t| > 0$ $\forall t$, which is sometimes called \textbf{weak survival}, ensures the infection persists globally, it does not require that the infection persists locally. However, we can consider the following stronger condition called \textbf{strong survival}. Suppose $x$ is some vertex (we will take it to be the origin on lattices) and let $\eta_t(x)$ be the state of vertex $x$ at time $t$. We say the process survives strongly if $\eta_t(x) = 1 \textrm{ u.o.}$ where $\textrm{u.o.}$ stands for ``unboundedly often," and we say an event $A_t$ occurs unboundedly often if $A_t$ occurs for an unbounded collection of times $t \in [0,\infty)$. We can then define the \textit{strong survival phase} as the set of $\lambda$ such that $\prob_{\lambda}(\eta_t(x) = 1 \textrm{ u.o.}) > 0$. Clearly strong survival implies weak survival, but the converse is not true in general \cite{Pemantle1992}. However, for the classical contact process on lattices $\mathbb{Z}^d, d \geq 1$, it is known that if the process survives weakly then it survives strongly \cite{LiggettSIS}.

One important property of the classical contact process is a type of monotonicity called \textbf{attractiveness} (or equivalently the process is called attractive.) Let $\eta^A_t$ denote the set of infected vertices at time $t$ starting from initial infected set $A \subseteq V$. Attractiveness states that given two initial infected vertex sets $A, B \subseteq V$ such that $A \subseteq B$, there exists a coupling such that $\eta_t^{A} \subseteq \eta^B_t$ for all $t \geq 0$. It is also true for the classical contact process that if $\lambda_1 \leq \lambda_2$, then for any initial infected set $A$

$$
\prob_{\lambda_1}(|\eta^A_t| > 0\textrm{ } \forall t) \leq \prob_{\lambda_2}(|\eta^A_t| > 0\textrm{ } \forall t).
$$
This allows us to characterize the existence of a phase transition with a \textbf{critical value} $\lambda_c$ such that for all $\lambda > \lambda_c, \prob_{\lambda}(|\eta^A_t| > 0$ $\forall t) > 0$ while for all $\lambda < \lambda_c, \prob_{\lambda_1}(|\eta^A_t| > 0$ $\forall t) = 0$. A phase transition is then said to exist if $0 < \lambda_c < \infty$. It is known that the classical contact process exhibits a phase transition on lattices $\mathbb{Z}^{d}, d \geq 1$, and that $\lambda_c$ scales like $1/d$ as $d \rightarrow \infty$. See \cite{LiggettSIS} for this and an overview of known results about the classical contact process on lattices. Attractivesness is an important ingredient in the the proofs of the existence of phase transitions and many other results about the classical contact process on lattices.

In contrast, the contact process with avoidance is not known to be attractive, which presents a notable challenge in studying it. The usual way to show the classical contact process is attractive is using the graphical construction we discuss in section \ref{sec:graphical}. \cite{CSW} show that this technique does not work for the contact process with avoidance; see their figure 2. Similarly, the contact process with avoidance is not known to exhibit ordering of survival probabilities in the infection parameter $\lambda$ for a fixed avoidance parameter $\alpha$. As a consequence, we cannot characterize the existence of a phase transition for the contact process with avoidance using a single critical value $\lambda_c$. Instead, \cite{CSW} show that on $\mathbb{Z}$, for every fixed $\alpha$ there exist values $0 < \lambda(\alpha)^- \leq \lambda(\alpha)^+ < \infty$ such that all $\lambda < \lambda(\alpha)^-$ are in the extinction phase and all $\lambda > \lambda(\alpha)^+$ are in the strong survival phase.

Our main results extend the study of the contact process with avoidance to $\mathbb{Z}^d, d \geq 2$. We show that for every fixed $\alpha > 0$ and $d \geq 2$, there exists $\lambda^{-}(\alpha,d) > 0$ and $\lambda^+(\alpha,d) < \infty$ such that on $\mathbb{Z}^d$ all $\lambda < \lambda(\alpha,d)^-$ are in the extinction phase and all $\lambda > \lambda(\alpha,d)^+$ are in the strong survival phase, and furthermore $\lambda(\alpha,d)^+$ scales like $1/d$ as $d \rightarrow \infty$. We also show that there exists a constant $C(\alpha)$ such that when $\lambda > C(\alpha)/d$ the process has a nontrivial invariant measure for $d$ sufficiently large. We conclude by discussing why our methods and most of our results also apply to the well-known SIRS model in which individuals become temporarily immune to reinfection after recovery and a related model called the contact process with isolation studied in \cite{CSW2}.

The contact process with avoidance, SIRS model, and contact process with isolation are natural extensions of the classical contact process. Indeed, interest in models with an isolated/immune state or dynamic edge behavior dates back to works such as \cite{DurrettNeuhauser} and \cite{GrossAdSIS}. Despite this, rigorous results about the survival and extinction of such models are still uncommon in the literature, in large part due to the technical challenges presented by these models' lack of attractiveness. We discuss notable results about these and related models in section \ref{sec:background}.

\subsection{Main Results}\label{sec:mainresults}

We formally define the \textbf{contact process with avoidance} $\{\xi_t\}_{t \geq 0}$ as follows. Let $G = (V,E)$ be a directed graph with vertices $V$ and edges $E$. Let $\xi_t(x) \in \{0,1\}$ denote the state of vertex $x$ at time $t$ where $0$ is healthy and $1$ is infected, and let $\xi_t(x,y)$ denote the state of directed edge $(x,y)$ where $1$ is active and $0$ is inactive. Let $\xi_t(V)$ denote the set of infected vertices at time $t$ and let $|\xi_t(V)|$ denote the size of this set. Given an initial configuration $\xi_0 \in \{0,1\}^V \times \{0,1\}^E$ and parameters $\lambda$ and $\alpha$, the process evolves according to the following rules:

\begin{enumerate}
    \item $\xi_t(x)$ goes from $0 \rightarrow 1$ at rate $\lambda \sum_{y \in V} \xi_t(y)\xi_t(y,x)\mathbf{1}_{(y,x) \in E}$.
    \item $\xi_t(x)$ goes from $1 \rightarrow 0$ at rate $1$.
    \item $\xi_t(x,y)$ goes from $1 \rightarrow 0$ at rate $\alpha$ if $\xi_t(x) = 1$ and $\xi_t(y) = 0$.
    \item $\xi_t(x,y)$ goes from $0 \rightarrow 1$ when $\xi_t(x) = 0$.
\end{enumerate}
Here rate means that the time to event follows an Exponential distribution with the given rate parameter. Typically, we will take the initial configuration $\xi_0$ to have infected set $V_0 = \{\mathbf{0}\}$ be the origin and initial active edge set $E_0 = E$ to allow the process to start with all edges active.

Because the contact process with avoidance is not known to be attractive, we cannot characterize its phase transitions with a single critical value $\lambda_c(d)$. Instead, we follow \cite{CSW} by defining upper and lower critical values for $G = \mathbb{Z}^d$ as follows.

\begin{equation}
    \begin{aligned}
    \lambda(\alpha,d)^- &= \inf\{\lambda: \prob_{\lambda}(|\xi_t(V)| > 0 \textrm{ } \forall t) > 0\},\\
    \lambda_w(\alpha,d)^+ &= \sup\{\lambda: \prob_{\lambda}(|\xi_t(V)| > 0 \textrm{ } \forall t) = 0\},\\
    \lambda(\alpha,d)^+ &= \sup\{\lambda: \prob_{\lambda}(\xi_t(\mathbf{0}) = 1 \textrm{ u.o.}) = 0\}.
    \end{aligned}
\end{equation}

Under these definitions, the process dies out almost surely for all $\lambda < \lambda(\alpha,d)^-$, survives weakly with positive probability for all $\lambda > \lambda(\alpha,d)^+$, and survives strongly with positive probability for all $\lambda > \lambda(\alpha,d)^+$ on $\mathbb{Z}^d$. However, we are unable to say precisely what happens for $\lambda(\alpha,d)^- \leq \lambda \leq \lambda_w(\alpha,d)^+$ and $\lambda_w(\alpha,d)^+ \leq \lambda \leq \lambda(\alpha,d)^+$. We now state our main results.

\begin{theorem}\label{thm:phase}
    Let $G = (V,E)$ where $V = \mathbb{Z}^d$ and $E = \{(x,y):||x-y||_1 = 1\}$ and fix $\alpha > 0$. Then
    \begin{enumerate}
        \item [(a)] $\lambda(\alpha,d)^- \geq \frac{1}{2d}$,
        \item [(b)] For all $d \in \mathbb{N}, \lambda(\alpha,d)^+ \leq a_1 + a_2 \alpha$ where $a_1$ and $a_2$ are the constants in \cite{CSW},
        \item [(c)] There exists a constant $D(\alpha) \in \mathbb{N}$ depending on $\alpha$ such that for all $d \geq D(\alpha)$, $\lambda(\alpha,d)^+ \leq C/d$ where the constant $C = 2((1 - e^{-1})(e^{-4\alpha})(e^{-3}))^{-1}$.
    \end{enumerate}
\end{theorem}

Theorem \ref{thm:phase} states that for any fixed $\alpha$, the contact process with avoidance has a phase transition $\lambda$ between almost sure extinction and positive probability of strong survival. Furthermore, the upper and lower critical values scale like $1/d$ as $d \rightarrow \infty$.

(a) and (b) are not difficult to prove. (a) follows from the fact that on any graph with maximum degree $M$, the classical contact process is subcritical for all $\lambda \leq 1/M$. This can be shown with a simple random walk comparison, see \cite{LiggettSIS}. The idea is that no matter the current configuration, when $\lambda \leq 1/M$ the next event that changes the size of the infected set is more likely to be a recovery than an infection and so the size of the infected set drifts toward zero. The same argument works for the contact process with avoidance, since the avoidance mechanism cannot lead to additional infections compared to the classical contact process. (b) follows from the results of \cite{CSW}, who note that their proof that the process can survive strongly on $\mathbb{Z}$ can also be applied to any graph that contains $\mathbb{Z}$ as a subgraph, which all $\mathbb{Z}^d, d \geq 1$ do.

Thus, the real work is in proving (c). Our primary tool is a graphical construction that allows us to compare the infection to oriented percolation in $\mathbb{Z}^d$, which is known to have a critical value that scales like $\frac{1}{d}$ \cite{CoxDurrett1983}. However, because the contact process with avoidance is not known to be attractive, our construction differs from that used to derive results about the classical contact process.

Another central result about the classical contact process is the \textbf{complete convergence theorem}. Let $\tau = \inf\{t \geq 0: |\eta_t| = 0\}$ be the time to extinction of the classical contact process $\eta_t$ on $\mathbb{Z}^d$. The complete convergence theorem states for any initial infected set $V_0$

\begin{equation}
\mu(\eta_t^{V_0}) \xrightarrow{w} \delta_0 \prob^{V_0}(\tau < \infty) + \bar{\nu} \prob^{V_0}(\tau = \infty)  
\end{equation}

where $\mu(\eta_t^{V_0})$ is the law of the process starting from initial infected set $V_0$, $\delta_0$ is point mass on $\nu_t = \emptyset$, $\bar{\nu}$ is the upper invariant measure obtained as the weak limit of the process starting from all vertices infected, and $P^{V_0}()$ denotes probability given initial infected set $V_0$. The proof of the complete convergence theorem relies heavily on the fact that the classical contact process is attractive and self-dual. Because the contact process with avoidance is not known to have these properties, we are unable to prove a fully analogous theorem. Instead, we prove the following weaker result.

\begin{theorem}\label{thm:invariant}   
    Fix $\alpha > 0, d \geq D(\alpha)$, and $\lambda > C/d$ where $C = 2\left((1-e^{-1})(e^{-4 \alpha}) (e^{-3}) \right)^{-1}$, and let $\{ \xi_t\}_{t \geq 0}$ be the contact process with avoidance on $\mathbb{Z}^d$ with initial configuration $\xi_0$ chosen so that all vertices are initially infected and all edges are initially active. Then $\xi_t$ has a non-trivial measure $\nu$ that is translation and time invariant.
\end{theorem}

Theorem \ref{thm:invariant} gives the existence of a non-trivial invariant measure $\nu$ in the regime where $d \geq D(\alpha
), \lambda > C/d$. Since the contact process with avoidance is not known to be attractive, our proof does not show that $\nu$ is unique and only shows that $\xi_t$ converges weakly to $\nu$ over a certain subsequence.

\subsection{Background and Related Work}\label{sec:background}

The classical contact process has been thoroughly studied on lattices $\mathbb{Z}^d, d \geq 1$. Again, we refer the reader to \cite{LiggettSIS}. In addition to the results about monotonicity, phase transitions, and the complete convergence theorem discussed in section \ref{sec:mainresults}, results exist about the growth of the infected set over time and behavior of the process at criticality. 

Most closely related to the contact process with avoidance are several other variant contact process models with dynamic edge behavior. One of the first such models proposed was the adaptive SIS or evoSIS model of \cite{GrossAdSIS} in which IS edges randomly rewire. This model has been the subject of substantial interest and studied, particularly in the physics literature, using mean field approximations, moment closures, and computational techniques \cite{GrossReview, GuoEtAl, PipatsartEtAl, DemirelEtAl}. However, it is has proven difficult to study rigorously and many fundamental questions remain open.

The contact process with avoidance, among other models, could be viewed as a modification of the evoSIS model for which rigorous results are more attainable, as the results of \cite{CSW} demonstrate. Another such model is the evoSIR model, which maintains dynamic edge rewiring but uses an underlying SIR model. In an SIR model vertices become permanently immune to infection upon recovery. Since each vertex can become infected at most once, SIR models tends to be simpler to study than SIS models. Accordingly, more is known about the evoSIR model, including the existence of phase transitions on Erd\"os-Renyi and configuration model random graphs \cite{JiangEtAl, DurrettYao}.

Another class of related models are those that decouple the dynamic edge behavior from the vertex states. The contact process with dynamic random edges, introduced in \cite{Linker2020}, features temporary edge avoidance similar to that the contact process with avoidance, but edges deactivate and reactive independently of the spread of infection. \cite{Linker2020} characterize phase transitions for this model on $\mathbb{Z}$, while \cite{Deshayes2026} extend the results of \cite{Linker2020} to $\mathbb{Z}^d$. Similarly, \cite{JacobMorters2017} and \cite{JacobMorters2025} study a variant of the evoSIS model on random graphs in which all edges rather than only IS edges randomly rewire. Unlike the evoSIS model and contact process with avoidance, these models tend to be attractive or at least stochastically ordered in some of their parameters.

Contact process variants not known to be attractive have also arisen from models that introduce additional vertex states. The well-known SIRS model in which infectious individuals become temporarily immune to reinfection upon recovery \cite{LamSIRS} and the contact process with isolation studied in \cite{CSW2} in which infected vertices can temporarily isolate until they recover are two such examples. \cite{DurrettNeuhauser} study the SIRS model on $\mathbb{Z}^2$ and show the existence of a phase transition and of a nontrivial stationary measure in the supercritical phase. However, their results rely on isoperimetric properties particular to $\mathbb{Z}^2$. \cite{GrimmettSIRS} build on an analogous result of \cite{Kuulasmaa} for the SIR model to show that for a fixed rate of transition from R to S, the SIRS model can survive indefinitely on $\mathbb{Z}^d$ when $d \geq 2$ if the infection rate is sufficiently large. However, their result gives weak survival, and so does not yield that the infection has the spatial recurrence properties needed to show the existence of a nontrivial invariant measure, nor do they quantitatively characterize value of the infection parameter necessary for survival. We discuss how our methodology and results apply to these models in section \ref{sec:isolation}.

As with dynamic edge behavior, attractiveness can be recovered by decoupling the isolated/immune state from the infection process. \cite{Remenik} studies a variant of the contact process in which all vertices can become temporarily isolated, finding it attractive in the usual sense and stochastically ordered in some of its parameters. He characterizes phase transitions on $\mathbb{Z}^d$ and derives a complete convergence theorem for this model.

\subsection{Graphical Construction}\label{sec:graphical}

A useful tool for studying the contact process and its variants is a graphical construction sometimes called the Harris construction in honor of Ted Harris. For the classical contact process, this construction is defined as follows. Given a graph $G = (V,E)$, consider the spacetime region $G \times [0,\infty)$. At each vertex $x$ and directed edge $e$ of the graph, we define a Poisson process of temporal marks as follows:

\begin{enumerate}
    \item[(C1)] On each directed edge $e\in G$, define a Poisson process on $\{e\}\times [0,\infty)$ with intensity $\lambda$ that generates infection arrows.
    \item[(C2)] On each vertex $x\in G$, define a Poisson process on $\{x\}\times [0,\infty)$ with intensity $1$ that generates recovery dots.
\end{enumerate}

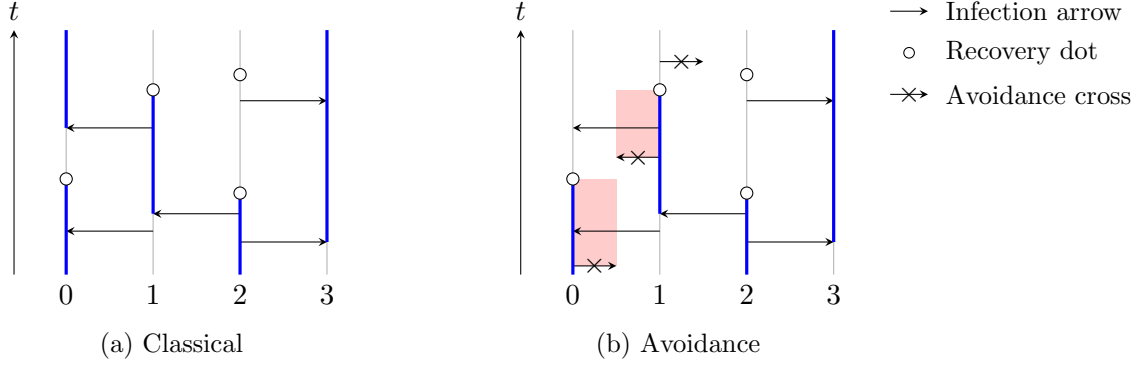
\begin{figure}[ht]
\centering

\begin{tabular}{@{}c@{\hspace{1.75cm}}c@{\hspace{0.30cm}}c@{}}

\begin{subfigure}[t]{0.30\linewidth}
\vspace{0pt}
\centering

\begin{tikzpicture}[scale=1.15, >=stealth]

    \def\T{2.8125}
    \def\rad{0.07}

    \foreach \x in {0,1,2,3}{
        \draw[gray!60, thin] (\x,0) -- (\x,\T);
        \node[below] at (\x,0) {$\x$};
    }

    \draw[->] (-0.6,0) -- (-0.6,\T)
        node[above] {$t$};


    \draw[->] (2,0.375) -- (3,0.375);

    \draw[->] (1,0.50) -- (0,0.50);

    \draw[->] (2,0.70) -- (1,0.70);

    \draw[->] (1,1.6875) -- (0,1.6875);

    \draw[->] (2,2.00) -- (3,2.00);


    \draw[blue, very thick]
        (0,0) -- (0,{1.10-\rad});

    \draw[blue, very thick]
        (0,1.6875) -- (0,\T);

    \draw[blue, very thick]
        (1,0.70) -- (1,{2.125-\rad});

    \draw[blue, very thick]
        (2,0) -- (2,{0.9375-\rad});

    \draw[blue, very thick]
        (3,0.375) -- (3,\T);

    \draw[fill=white] (0,1.10) circle (\rad);
    \draw[fill=white] (1,2.125) circle (\rad);
    \draw[fill=white] (2,0.9375) circle (\rad);

    \draw[fill=white] (2,2.30) circle (\rad);

\end{tikzpicture}

\caption{Classical}
\label{fig:1a}
\end{subfigure}

&

\begin{subfigure}[t]{0.30\linewidth}
\vspace{0pt}
\centering

\begin{tikzpicture}[scale=1.15, >=stealth]

    \def\T{2.8125}
    \def\rad{0.07}

    \foreach \x in {0,1,2,3}{
        \draw[gray!60, thin] (\x,0) -- (\x,\T);
        \node[below] at (\x,0) {$\x$};
    }

    \draw[->] (-0.6,0) -- (-0.6,\T)
        node[above] {$t$};

    \fill[red, fill opacity=0.20]
    (0,0.1) rectangle (0.5,1.10);

    \fill[red, fill opacity=0.20]
    (0.5,1.35) rectangle (1,2.125);


    \draw[->] (2,0.375) -- (3,0.375);

    \draw[->] (1,0.50) -- (0,0.50);

    \draw[->] (2,0.70) -- (1,0.70);

    \draw[->] (1,1.6875) -- (0,1.6875);

    \draw[->] (2,2.00) -- (3,2.00);


    \draw[->] (0,0.10) -- (0.50,0.10)
        node[midway] {$\times$};

    \draw[->] (1,1.35) -- (0.50,1.35)
        node[midway] {$\times$};

    \draw[->] (1,2.45) -- (1.50,2.45)
        node[midway] {$\times$};


    \draw[blue, very thick]
        (0,0) -- (0,{1.10-\rad});

    \draw[blue, very thick]
        (1,0.70) -- (1,{2.125-\rad});

    \draw[blue, very thick]
        (2,0) -- (2,{0.9375-\rad});

    \draw[blue, very thick]
        (3,0.375) -- (3,\T);

    \draw[fill=white] (0,1.10) circle (\rad);
    \draw[fill=white] (1,2.125) circle (\rad);
    \draw[fill=white] (2,0.9375) circle (\rad);

    \draw[fill=white] (2,2.30) circle (\rad);

\end{tikzpicture}

\caption{Avoidance}
\label{fig:1b}
\end{subfigure}

&

\begin{minipage}[t]{0.20\linewidth}
\vspace{0pt}

\begin{tikzpicture}[>=stealth]

    \draw[->] (0,0) -- (0.50,0);
    \node[right] at (0.60,0)
        {\small Infection arrow};

    \draw[fill=white] (0.25,-0.55) circle (0.07);
    \node[right] at (0.60,-0.55)
        {\small Recovery dot};

    \draw[->] (0,-1.10) -- (0.50,-1.10)
        node[midway] {$\times$};
    \node[right] at (0.60,-1.10)
        {\small Avoidance cross};

\end{tikzpicture}

\end{minipage}

\end{tabular}

\caption{Examples of realized graphical constructions for the classical contact process and contact process with avoidance.}

\end{figure}

We can then realize the process on $G \times [0,\infty)$ given some initial state $\xi_0\in \{0,1\}^V$ using these marks by doing the following. See figure \ref{fig:1a} for an example realization.

\begin{enumerate}
    \item Label infected vertices (1, blue) vertically in time until a recovery dot is reached.
    \item When an infection arrow is observed, if the source vertex is infected (1, blue), infect the target vertex and repeat step 1. for this newly infected vertex.
\end{enumerate}

To define a similar construction for the contact process with avoidance, we add the following to C1 and C2

\begin{itemize}
    \item[(A3)] On each directed edge $e\in G$, define a Poisson process on $\{e\}\times [0,\infty)$ with intensity $\alpha$ that generates avoidance crosses.
\end{itemize}

We can then realize the process on $G \times [0,\infty)$ given some initial state $\xi_0\in \{0,1\}^V \times \{0,1\}^E$ using these marks by doing the following. See figure \ref{fig:1b} for an example realization.

\begin{enumerate}
    \item Label infected vertices (1, blue) vertically in time until a recovery dot is reached.
    \item Whenever a recovery dot is observed on a vertex $x$, set all directed edges $(x,y)$ where $y \sim x$ to state ($1$, uncolored).
    \item When an avoidance cross is observed on a directed edge $(x,y)$, if vertex $x$ is in state ($1$, blue) and vertex $y$ is in state $0$, set edge $(x,y)$ to state $0$ and color it red until a recovery dot is observed on vertex $x$.
    \item When an infection arrow is observed, if the source vertex is infected (1, blue) and the directed edge is not (0, red), infect the target vertex and repeat steps 1. and 2. for this newly infected vertex.
\end{enumerate}

\section{Survival of the Contact Process with Avoidance on $\mathbb{Z}^d$}\label{sec:block}
In this section we prove (c) of Theorem \ref{thm:phase}. As with the classical contact process, our goal is to define a block construction in which we can compare spacetime regions to sites in a oriented site percolation. However, the lack of attractiveness makes the construction of appropriate regions more complicated.

For the contact process with avoidance on the lattice $\mathbb{Z}^d$, define the region $R(j,t)$ to consist of a center vertex $j$ and all the directed edges $(i,j)$ from each of the $2d$ neighbors $i$ of $j$ to $j$. We define the temporal center of the region at time $k \Delta t$ take increments of $\Delta t$ that we will apply in different ways to the vertices and edges. We include in $R(j,t)$
\begin{enumerate}
    \item The Poisson process of infection arrows on each edge $(i,j)$ from time $k \Delta t +\Delta t, k \Delta t + 2\Delta t$
    \item The Poisson process of avoidance crosses on each edge $(i,j)$ from time $k \Delta t -2\Delta t, k \Delta t + 2\Delta t$
    \item The Poisson process of recovery dots on $j$ from $k \Delta t, k \Delta t + 4\Delta t$
\end{enumerate}

Suppose we start from the configuration where the origin $v_O$ is infected, all other vertices are healthy, and no edges are avoiding. We first treat the region $R(v_O,0)$ as a special case and say it is \textbf{good} if there are no recovery dots on $v_O$ during the time interval $(0,4 \Delta t)$. For each $k \in \{1,2,3\ldots\}$ we say the region $R(j,2k\Delta t)$ is \textbf{good} if all of the following:
\begin{itemize}
    \item[N1] There is some vertex $i$ that is a neighbor of $j$ such that $R(i,2(k-1) \Delta t)$ is good and there is an infection arrow on $(i,j)$ in the time interval $(2k+1)\Delta t, (2k + 2)\Delta t$.
    \item[N2] There is recovery dot on $j$ in the time interval $2k \Delta t, (2k+ 1)\Delta t$
    \item[N3] There are no avoidance crosses on $(i,j)$ in the time interval $(2k - 2)\Delta t, (2k + 2)\Delta t$
    \item[N4] There are no recovery dots on $j$ in the time interval $(2k + 1)\Delta t, (2k + 4)\Delta t$
\end{itemize}

To understand the utility of this construction, suppose $R(v_O,0)$ is good and there is some neighbor of the origin $j$ satisfying N1-N4 for $k = 1$ with the infection arrow in N1 coming from $v_O$ so that $R(j,2\Delta t)$ is good. N1, N3, and N4 ensure that $v_O$ successfully infects $j$ between time $3 \Delta t$ and $4 \Delta t$ without being blocked by an avoidance on the edge $(v_O,j)$ and that $j$ remains infected until time $6 \Delta t$. N2 ensures that when we consider the region $R(\ell,4\Delta t)$ for each neighbor $\ell$ of $j$, the edge $(j,\ell)$ cannot be avoiding due to any avoidance crosses that originated before time $2 \Delta t$. If we then consider some neighbor $\ell$ of vertex $j$ satisfying N1-N4 for $k = 2$ so that $R(\ell,4\Delta)$ is good, we see these properties propogate. N1, N3, and N4 ensure that $j$ successfully infects $\ell$ between time $5 \Delta t$ and $6 \Delta t$ without being blocked by an avoidance on the edge $(j,\ell)$ and that $\ell$ remains infected until time $8 \Delta t$. N2 ensures that when we consider the region $R(h,6\Delta t)$ for each neighbor $h$ of $\ell$, the edge $(\ell,h)$ cannot be avoiding due to any avoidance crosses that originated before time $4 \Delta t$. We can continue in this manner, alternating between regions centered at even vertices when $k$ is even and regions centered at odd vertices when $k$ is odd, noting that good regions allow the infection to propagate forward through time.

We now consider the probability that a region is good. $R(v_O,0)$ is good if there are no recovery dots on $v_O$ during the time interval $(0,4 \Delta t)$ which occurs with probability $p_0 := e^{-4\Delta t}$.

For $k \geq 1$ the probability that a region $R(j,2k\Delta t)$ depends on if there is at least one neighbor $i$ of $j$ such that the region $R(i,2(k-1)\Delta t)$ is good. If no such $i$ exists, then $R(j,2k\Delta t)$ is good with probability $0$. Otherwise, if at least one such $i$ exists, $R(j,2k\Delta t)$ is good with probability at least

\begin{equation}\label{eq:goodprob1}
p_G := (1- e^{-\lambda\Delta t})(1 - e^{-\Delta t})(e^{-4\alpha \Delta t})(e^{-3 \Delta t}).  
\end{equation}

Suppose we fix $\alpha > 0$, let $\Delta t = 1$, and let $\lambda = C/d$ for a constant $C$ to be chosen later. Using Taylor expansion, as $d \rightarrow \infty$

\begin{equation}\label{eq:goodprob2}
p_G = (1- e^{-C/d})(1 - e^{-1})(e^{-4\alpha})(e^{-3}) = \frac{C((1 - e^{-1})(e^{-4\alpha})(e^{-3}))}{d} + O(d^{-2}).
\end{equation}

Thus if we choose $C = 2((1 - e^{-1})(e^{-4\alpha})(e^{-3}))^{-1}$, then there exists $D_1(\alpha)$ such that for all $d \geq D_1(\alpha)$, $p_G > \frac{3}{2d}$.

We now compare our path of good regions to the percolation cluster in an oriented site percolation in dimension $\mathbb{Z}^{d+1}$ where $d > D_1(\alpha)$. Recall that in this model, the $+1$ dimension is discrete time, and we consider sites on the non-negative orthant of $\mathbb{Z}^d$. Each site is occupied independently with probability $p$; here we will take $p = p_G$. Since $p_G > \frac{3}{2d}$, there exists $D_2$ such that for all $d \geq D_2$ the percolation is supercritical and there is positive probability that there exists an infinite path in time of occupied sites starting from the origin.

We map our spacetime regions onto sites in this oriented percolation. We would like to exhibit a coupling where good regions dominate the occupied sites in the percolation model. However, since a region can only be good if it neighbors another good region in the previous time increment, we require some modifications first. Suppose we realize the percolation model in the usual way and then do the following.
\begin{enumerate}
    \item Remove occupancy from all sites except the origin at time $0$.
    \item Proceding forward in discrete time units, remove occupancy from any site that does not have at least one occupied neighbor at the previous time increment.
\end{enumerate}

We first observe that if there is an infinite path in time of occupied sites starting from the origin in the original percolation model, it will remain after this modification, since the origin must have been occupied at time $0$ and any site in the path must have at least one occupied neighbor in the previous time increment. In addition, we can now couple our graphical construction with the modified oriented percolation model so that the good regions in the graphical construction stochastically dominate the occupied sites in the modified oriented site percolation model. And since by our choice of $p_G$, the modified oriented site percolation has positive probability of an infinite path in time of occupied sites starting from the origin, there is positive probability that the infection in the contact process with avoidance survives forever. That is, for all $d \geq D(\alpha) := \max\{D_1(\alpha),D_2\}$, $\prob_{\lambda}(|\xi_t^{v_O}| > 0 \textrm{ } \forall t ) > 0$ for all $\lambda > C/d$.

To obtain strong survival, we apply the results of \cite{bezuidenhout1990critical} to our block construction. They show that in a supercritical oriented site percolation, if there is an infinite path from the origin then the sites corresponding to the regions $R(v_O,2k \Delta t)$ are occupied for infinitely many $k \in \mathbb{N}$ and thus $\prob^{\lambda}(\xi_t(\mathbf{0}) = 1 \textrm{ u.o.}) > 0$.

\section{Existence of a Nontrivial Invariant Measure}

We now prove Theorem \ref{thm:invariant}. Fix $\alpha > 0, d \geq D(\alpha)$ and $\lambda > C/d$ where $C = 2((1 - e^{-1})(e^{-4\alpha})(e^{-3}))^{-1}$. Our construction comes as a corollary of the comparison to oriented percolation used in the proof of Theorem \ref{thm:phase}.

Recall in our notation $\{ \xi_t \}_{t \geq 0}$ is the contact process with avoidance on the configuration space,
\begin{equation}
    \Sigma = \left\{ \{ 0, 1\} \times \{0,1 \}^{2d} \right\}^{\mathbb{Z}}
\end{equation}
where $\xi_t(x)$ is the state of vertex $x$ at time $t$ with $\xi_t(x) = 1$ if $x$ is infected and $0$ is $x$ is healthy and $\xi_t(x,y)$ is the state of edge $(x,y)$ at time $t$ with $\xi_t(x,y) = 1$ if $(x,y)$ is active and $0$ is $(x,y)$ is inactive. We take the initial configuration $\xi_0$ to have all vertices infected and all edges active. 

Let $B_n$ to be ball of radius $n$ in the $L^\infty$ norm around the origin. For $n \in \{1,2,\ldots, \}$ we define the probability measure $\nu^{(n)}_t$ to be the joint distribution of the finite collection of variables,
\begin{equation}
    \left\{ \xi_t(x) \mid x \in B_n \right\} \cup \left\{  \xi_t(x,y) \mid x \in B_n, \ y \sim x\right\}.
\end{equation}
That is, $\nu^{(n)}_t$ is the marginal distribution of the statuses of all vertices and directed edges outgoing from vertices in $B_n$.

Our strategy is to construct a consistent sequence of marginal distributions of the contact process with avoidance on the balls $B_n$. Each collection $\{ \nu^{(n)}_t \}_{t \geq 0}$ has a bounded support across all $t \geq 0$ and therefore is tight. We then use the Kolmogorov extension theorem to take the natural extension of this sequence to $\Sigma$. In order to show our measure is time-invariant we consider the process run for an amount of time uniformly distributed over interval $[0,k]$. This is motivated by the fact that the sum of a Uniform$[0,k]$ random variable and a constant $t>0$ is nearly equal in distribution to a Uniform$[0,k]$ random variable when $k$ is large. Thus, running the contact process with avoidance for additional time $t>0$ after running it for a long, uniformly distributed time has only a small impact on its law.

Let $U_k, k \in \mathbb{N}$ be independent uniform variables supported on the intervals $[0,k]$. Starting from the ball $B_1$, we can use Prokhorov's theorem to form a sequence $\{ a^{(1)}_k \}_{k=1}^{\infty}$ with $a^{(1)}_k \rightarrow \infty$ so that the sequence of measures,
\begin{equation}
\label{nu1}
\left\{ \nu^{(1)}_{U_{a^{(1)}_k}} \right\}_{k=1}^\infty
\end{equation}
is convergent in distribution to the measure $\nu^{(1)}$ of infection and edge statuses on $B_1$. We can apply Prokhorov's theorem again to form a subsequence $\{ a^{(2)}_k \}_{k=1}^\infty$ of $\{ a^{(1)}_k \}_{k=1}^\infty$ so that,
\begin{equation}
\left\{ \nu^{(2)}_{U_{a^{(2)}_k}} \right\}_{k=1}^\infty
\end{equation}
is convergent in distribution to a measure $\nu^{(2)}$ supported on the infection and edge status of sites in the ball $B_2$. Since $\{a_k^{(2)} \}$ is a subsequence of $\{ a^{(1)}_k\}$, we conclude $\nu^{(2)}$ and $\nu^{(1)}$ agree on $B_1$. We now construct the sequence of measures $\nu^{(n)}$ consisting of the joint distribution of infection and edge statuses within the ball $B_n$ as follows. If we have a sequence $\{ a^{(n)}_k \}_{k =1}^\infty$ so that $a^{(n)}_k \rightarrow \infty$ and the collection $\{ \nu^{(m)}\}_{m=1}^{n}$ is a consistent sequence of distributions supported on the finite dimensional distributions on balls $B_m$ for $m = 1, \ldots, n$, we use Prokhorov's theorem to find a subsequence $\{ a^{(n+1)}_k\}_{k=1}^{\infty}$ so that $\left\{ \nu^{(n+1)}_{U{a^{(n+1)}_k}} \right\}_{k=1}^\infty$ is convergent in distribution to a measure $ \nu^{(n+1)} $ on finite dimensional distributions over $B_{n+1}$, and the collection of measure $\{ \nu^{(m)}\}_{m=1}^{n+1}$ is a consistent sequence of distributions. We note that for each $n$ the measure $\nu^{(n)}$ is supported on the set,
\begin{equation}
    \Sigma^{(n)} = \left\{ \{0, 1\} \times \{0,1\}^{2d} \right\}^{(2n+1)^d}.
\end{equation}

Applying this process iteratively we obtain the consistent sequence of distributions $\{ \nu^{(n)}\}_{n=1}^\infty $. Using the Kolmorogov extension theorem we can form a unique probability measure $\nu$ supported the space $\Sigma$ equipped with the smallest $\sigma$-field generated by the joint distribution of the infection and edge statuses of sites over balls $B_n$ for $n = \{1,2,\ldots \}$ which we denote by $\mathcal{F}$.

We now show that $\nu$ is non-trivial, translation invariant, and time invariant. Since our initial configuration $\xi_0$ is the configuration with all sites infected and no edges in avoidance, we can leverage the percolation comparison process from Theorem \ref{thm:phase} along with results from \cite{bezuidenhout1990critical} applied to oriented site percolation which when translated to the notation of our construction yield that
\begin{equation}
\liminf_{k \rightarrow \infty} \prob(R(v_O,2k \Delta t) \textrm{ is good}) > 0,
\end{equation}
and thus
\begin{equation}
    \label{liminfprob}
    \liminf_{t \rightarrow \infty} \mathbb{P}\left( \xi_t(\mathbf{0}) = 1 \right) >0.
\end{equation}

Letting $A = \{ \xi \in \Sigma \mid \xi(\mathbf{0}) = 1 \}$, we apply our construction and (\ref{liminfprob}) to conclude,
\begin{equation}
    \nu(A) = \lim_{k \rightarrow \infty} \mathbb{P}{(\xi_{a^{(1)}_k}}(\mathbf{0}) = 1) > 0,
\end{equation}
and thus the measure $\nu$ is nontrivial. By symmetry of the starting configuration $\xi_0$, we can conclude the law of the configuration $\xi_t$ is translation invariant for all $t > 0$. To see that $\nu$ is translation invariant, let $C$ be any event that is fully determined by the infection and avoidance statuses within the ball $B_n$, and let $x \in \mathbb{Z}^d$ with $\lVert x-y \rVert_\infty = m$. Let $C^x$ be the event $C$ applied to the process translated by $x$. It is clear by translation invariance of $\xi_t$ for all $t$ that,
\begin{equation}
    \nu^{(m+n)}(C) = \nu^{(m+n)}(C^x).
\end{equation}
Since $\mathcal{F}$ is the product $\sigma$ field generated by finite dimensional distributions over balls $B_n$, we conclude $\nu$ is translation invariant. 

Lastly, we show that $\nu$ is time invariant. We fix $t>0$, let $\tilde{\xi}_t$ be the contact process with avoidance with $\tilde\xi^{}_0$ sampled via $\nu$, and let $\tilde{\nu}_t$ be the law of $\tilde{\xi}^{}_t$. We will now show $\nu = \tilde{\nu}_t$ for all $t > 0$.

Fix $t > 0$. Let $D$ be any finite cylinder event that is fully determined the statuses of the vertices and edges in $B_n$. Let $\varepsilon > 0$ be a small and positive constant so that $t > \varepsilon$, and let $m = m(n, \varepsilon)$ be a constant depending on $\varepsilon$ and $n$ that we will fully specify later. From the definition of the sequence $a_k^{(m)}$, we can choose a constant $k = k(m, \varepsilon)$ large enough so that $a^{(m)}_k > \varepsilon^{-1}$ and sufficiently large to allow us form a joint distribution of
\begin{equation}
    \label{jointconstruction}
    \left(\xi_{U_{a^{(m)}_k}}, \tilde{\xi}_0 \right)
\end{equation}
so that the event $E = \{\xi_{U_{a^{(m)}_k}} \textrm{ and }\tilde{\xi}_0 \textrm{ agree on the ball }B_m\}$ occurs with probability at least $1-\varepsilon$. 

Next, we evolve $\{ \xi_{U_{a^{(m)}_k}+s}\}_{0 \leq s \leq t}$ and $\{\tilde{\xi_s}\}_{0 \leq s \leq t}$ using the same set of marks in the graphical construction. Given the event $E$, we construct a coupled region of sites between $\{ \xi_{U_{a^{(m)}_k}+s}\}_{0 \leq s \leq t}$ and $\{\tilde{\xi_s}\}_{0 \leq s \leq t}$ that we denote as $\{ W_s \}_{0 \leq s \leq t} \subseteq \mathbb{Z}^d$ as follows. We set the initial coupled region to be $W_0 =B_m$, where both processes agree on the event $E$. Let $\partial W_s$ denote the vertex boundary of $W_s$. We now use the following rules to evolve $\{ W_s \}_{0 \leq s \leq t}$,
\begin{enumerate}
    \item If $x \in \partial W_{s^-}$ and there is a neighboring vertex $y \notin W_{s^-}$ such that there is an infection arrow on the edge $(y,x)$ at time $s$, we remove $x$ and all edges pointing outward from $x$ from $W_s$.
    \item If $x \in \partial W_{s^-}$ and there is a neighboring site $y \notin W_{s^-}$ such that there is an avoidance cross on the edge $(y,x)$ at time $s$, we remove $x$ and all edges pointing outward from $x$ from $W_s$.
\end{enumerate}
The rules above guarantee that for any $0 \leq s \leq t$ and $x \in W_s$, the infection status of $x$ is shared between $\xi_{U_{a^{(m)}_k}+s}$ and $\tilde{\xi_s}$ and that for any $x \in W_s$, and $y \sim x$ that the avoidance status between the directed edge $(x,y)$ is shared between $\xi_{U_{a^{(m)}_k}+s}$ and $\tilde{\xi_s}$. Note that each vertex $x \in \partial W_s$ leaves the coupled region with rate,
\begin{equation}
    \sum_{y \sim x} (\lambda + \alpha) \mathbbm{1} \{ y \notin W_s \}.
\end{equation}
We can therefore compare the decay of our coupled region $\{ W_s \}_{0 \leq s \leq t}$ to Richardson's growth model \cite{Richardson1973} with growth rate $(\lambda + \alpha)$, which on $\mathbb{Z}^d$ is equivalent the classical contact process without recovery and with infection rate $\lambda + \alpha$. In this model, each healthy vertex $x$ becomes infected at rate $\sum_{y \in V}(\lambda + \alpha)\mathbb{1}_{(y,x) \in E}$ and subsequently remains infected forever. 

Let $F$ be the event that the ball $B_n$ is contained within $W_t$. Conditioned on $E$, we see $F^c$ occurs if and only if there exists a dual path of Richardson's growth model beginning in $B_n$ that reaches the set $B_m$ by time $t$. We can thus apply estimates from Theorem 1 of \cite{durrett1982contact} to choose $m = m(n, \varepsilon)$ sufficiently large so that
\begin{equation}
    \mathbb{P}(F \mid E) > 1- \varepsilon.
\end{equation}
Since on the event $F$ we have that $\xi_{U_{a^{(m)}_k}+t}$ and $\tilde{\xi_t}$ are coupled within $B_n$ we can conclude,
\begin{equation}
    \label{trianglesetup1}
    | \nu^{(m)}_{U_{a_k}^{(m)}+t}(D) - \tilde{\nu}_t(D) | < 2\varepsilon.
\end{equation}
We now compare the probabilities $\nu^{(m)}_{U_{a_k}^{(m)}+t}(D)$ and $\nu^{(m)}_{U_{a_k^{(m)}}}(D)$. Using that $U_{a_k^{(m)}}$ is uniformly distributed over the interval $[0, a^{(m)}_k]$, we can apply the tower property to compute,
\begin{equation}
    \label{towerproperty1}
    \nu^{(m)}_{U^{(m)}_{a_k}}(D) = \int_{0}^{a^{(m)}_k} \frac{1}{a^{(m)}_k} \nu^{(m)}_s(D) \ ds.
\end{equation}
We can similarly write the probability,
\begin{equation}
    \label{towerproperty2}
\nu^{(m)}_{U_{a^{(m)}_k}+t}(D) = \int_{t}^{a^{(m)}_k+t} \frac{1}{a^{(m)}_k} \nu^{(m)}_s (D) \ ds.
\end{equation}
We can then combine (\ref{towerproperty1}), (\ref{towerproperty2}), and the triangle inequality to bound
\begin{equation}
    \label{trianglesetup2}
    \begin{aligned}
        &\left|\nu^{(m)}_{U^{(m)}_{a_k}}(D)  - \nu^{(m)}_{U_{a^{(m)}_k}+t}(D)   \right| \\
        & \leq 0 +  \int_0^t \left| \frac{1}{a^{(m)}_k} (1) \right| ds+ \int_{a^{(m)}_k}^{a^{(m)}_k+t}  \left| \frac{1}{a^{(m)}_k} (1) \right|ds \\
        & \leq 2 \varepsilon t,
    \end{aligned}
\end{equation}
with the last inequality holding using the bound $a^{(m)}_k > \varepsilon^{-1}$.

Lastly we have defined $k = (m,\varepsilon)$ to be sufficiently large so that,
\begin{equation}
    \label{trianglesetup3}
    \left| \nu(D) -  \nu^{(m)}_{U_{a^{(m)}_k }}(D)  \right| < \varepsilon.
\end{equation}
Using (\ref{trianglesetup1}), (\ref{trianglesetup2}), (\ref{trianglesetup3}) and the triangle inequality we have,
\begin{equation}
    | \nu(D) - \tilde{\nu}_t(D) | < (3+2t) \varepsilon.
\end{equation}
Since $\varepsilon > 0$ is an arbitrarily small constant we conclude $\nu(D) = \tilde{\nu}_t(D)$. Lastly, since $\nu$ and $\tilde{\nu}_t$ are constructed using the smallest $\sigma$ field generated by all finite cylinder events, we conclude $\nu = \tilde{\nu}_t$, and thus $\nu$ is time invariant. 

\section{Models with Vertex Isolation}\label{sec:isolation}

Most of our results also hold for the contact process with isolation studied by \cite{CSW2} and the SIRS model, namely (a) and (c) of Theorem \ref{thm:phase} for a different constant $C$ and dependence of $D(\alpha)$ on the appropriate parameter for each model and Theorem \ref{thm:invariant}. For the contact process with avoidance $\alpha$ controls the rate at which infected vertices become isolated, while for the SIRS model $\alpha$ controls the rate at which vertices transition from R to S. Although it is possible to make some statement about weak survival for $d \geq 2$ as \cite{GrimmettSIRS} do for the SIRS model, we are unable to claim (b) of Theorem \ref{thm:phase} for these models because their behavior on $\mathbb{Z}$ is not fully understood. \cite{Heeszel2026} showed that if the SIR model survives weakly on $\mathbb{Z}$ then it also survives strongly, and the techniques of \cite{CSW} can be used to show the SIRS model survive strongly when $\alpha$ is sufficiently large. However, fully characterizing the survival and extinction of the contact process with isolation and the SIRS model on $\mathbb{Z}$ remains an open problem.  The methodology for proving (a) and (c) in both cases is analogous to what we used for the contact process with avoidance. To avoid repeating the same arguments with minor variations, we sketch some details only for the contact process with isolation. 

Formally, the contact process with isolation $\{\zeta_t\}_{t \geq 0}$ is defined as follows. Let $G = (V,E)$ be a graph with vertices $V$ and edges $E$ and let $\zeta_t(x) \in \{0,1,-1\}$ denote the state of vertex $x$ at time $t$ where $0$ is healthy, $1$ is infected, and $-1$ is isolated. Given initial configuration $\zeta_0 \in \{0,1,-1\}^V$ and parameters $\lambda$ and $\alpha$, the process evolves according to the following rules.

\begin{enumerate}
    \item $\zeta_t(x)$ goes from  $0 \rightarrow 1$ at rate $\lambda \sum_{y \in V}\mathbf{1}_{\zeta_t(y) = 1} \mathbf{1}_{(y,x) \in E}$,
    \item $\zeta_t(x)$ goes from $1 \rightarrow -1$ at rate $\alpha$,
    \item $\zeta_t(x)$ goes from $1,-1 \rightarrow 0$ at rate $1$.
\end{enumerate}

\noindent Again rate means that the time to event follows an Exponential distribution with the given rate parameter. In this model, the parameter $\alpha$ controls the isolation rate rather than the avoidance rate.

(a) of Theorem \ref{thm:phase} holds for the contact process with isolation by the same argument given in section \ref{sec:mainresults} for the contact process with avoidance. (c) of Theorem \ref{thm:phase} and Theorem \ref{thm:invariant} hold because we can create a block construction analogous to the one used in section \ref{sec:block}. We first define a graphical construction for the contact process with isolation on spacetime region $G \times [0,\infty)$ by combining rules C1 and C2 from section \ref{sec:graphical} with the following rule

\begin{enumerate}
    \item[(I3)] On each vertex $x \in G$, define a Poisson process on $x \times \{0,\infty\}$ with intensity $\alpha$ that generates isolation crosses.
\end{enumerate}

We can then realize the the contact process with isolation on $G \times [0,\infty)$ given some initial configuration $\zeta_0\in \{0,1,-1\}^V$ using these marks by doing the following.

\begin{enumerate}
    \item Label infected vertices $1$ blue vertically in time until a recovery dot or isolation cross is reached.
    \item Whenever a recovery dot is observed on a vertex $x$, set the state of $x$ to (0, uncolored).
    \item When an isolation cross is observed on a vertex $x$, if $x$ is in state (1,blue), set $x$ to state (-1,red) and color $x$ red vertically in time until a recovery dot is reached.
    \item When an infection arrow is observed, if the source vertex is infected (1, blue) and the target vertex is healthy (0, uncolored), infect the target vertex and repeat steps 1., 2., and 3. for this newly infected vertex.
\end{enumerate}

See \cite{CSW2} for a more detailed description of this construction. We can then define a block construction analogous to the one in \ref{sec:block}. We define $R(j,t)$ in the same manner except that we include

\begin{enumerate}
    \item The Poisson process of infection arrows on each edge $(i,j)$ from time $k \Delta t +\Delta t, k \Delta t + 2\Delta t$
    \item The Poisson process of isolation crosses on $j$ from time $k \Delta t, k \Delta t + 4\Delta t$
    \item The Poisson process of recovery dots on $j$ from $k \Delta t, k \Delta t + 4\Delta t$
\end{enumerate}

We say the region $R(v_O,0)$ is good if there are no recovery does or isolation crosses on $v_O$ in the time interval $(0,4 \Delta t)$. For each $k \in \{1,2,3\ldots\}$ we say the region $R(j,2k\Delta t)$ is \textbf{good} if and all of the following:
\begin{itemize}
    \item[M1] There is some vertex $i$ that is a neighbor of $j$ such that $R(i,2(k-1) \Delta t)$ is good and there is an infection arrow on $(i,j)$ in the time interval $(2k+1)\Delta t, (2k + 2)\Delta t$.
    \item[M2] There is recovery dot on $j$ in the time interval $2k \Delta t, (2k+ 1)\Delta t$
    \item[M3] There are no isolation crosses on $j$ in the time interval $2k\Delta t, (2k + 4)\Delta t$
    \item[M4] There are no recovery dots on the vertex $j$ in the time interval $(2k + 1)\Delta t, (2k + 4)\Delta t$
\end{itemize}

Analogous to before, M1-M4 ensure that $j$ successfully receives the infection from some neighbor $i$ and is in position to spread the infection to regions centered on its neighbors in the next time step. We then proceed as in section \ref{sec:block} with this definition of a good region with some small differences in the calculations in equations \ref{eq:goodprob1} and \ref{eq:goodprob2} and thus in the values of $C$ and $D(\alpha)$. Since Theorem \ref{thm:invariant} follows from the block construction, it also holds for the contact process with isolation.

It is possible to go through the same procedure again for the SIRS model. We first formally define the SIRS model $\{\Upsilon_t\}_{t \geq 0}$ as follows. Let $G = (V,E)$ be a graph with vertices $V$ and edges $E$ and let $\Upsilon_t(x) \in \{0,1,-1\}$ denote the state of vertex $x$ at time $t$ where $0$ is healthy, $1$ is infected, and $-1$ is immune (isolated.) Given an initial configuration $\Upsilon_0 \in \{0,1-1\}^V$ and parameters $\lambda$ and $\alpha$, the process evolves according the following rules.

\begin{enumerate}
    \item $\Upsilon_t(x)$ goes from $0 \rightarrow 1$ at rate $\lambda \sum_{y \ in V}\mathbf{1}_{\Upsilon_t(y) = 1}\mathbf{1}_{(yx,) \in E}$,
    \item $\Upsilon_t(x)$ goes from $1 \rightarrow -1$ a rate $\alpha$,
    \item $\Upsilon_t(x)$ goes from $-1 \rightarrow 0$ at rate $1$.
\end{enumerate}

\noindent Essentially, the only difference between the contact process with isolation and the SIRS model is whether infected vertices can return to the healthy state without first becoming immune (isolated) for some time.

To show (a) and (c) of Theorem \ref{thm:phase} and Theorem \ref{thm:invariant} for the SIRS model, we modify the graphical construction and definition of a good region to use marks that control transitions from $1$ to $-1$ in place of avoidance crosses and marks that control transitions from $-1$ to $0$ in place of recovery dots. To avoid repeating ourselves, we omit the details. We refer to the reader to \cite{CSW2} for further discussion of the similarities between the SIRS model and contact process with isolation.

\bibliographystyle{abbrv}
\bibliography{ref}

\end{document}